\documentclass[letterpaper, 10 pt, conference]{ieeeconf}

\usepackage{cite}
\usepackage{amsmath,amssymb,amsfonts}
\usepackage{algorithmic}
\usepackage{graphicx}
\usepackage{subcaption}
\usepackage{textcomp}
\usepackage{algorithm}
\usepackage{algorithmic}
\usepackage{lipsum}
\usepackage{dblfloatfix}

\usepackage{graphicx}
\usepackage{epsfig} % for postscript graphics files
\usepackage{cite}
\usepackage{booktabs}
\usepackage[dvipsnames]{xcolor}

\IEEEoverridecommandlockouts                              % This command is only needed if you want to use the \thanks command

\newtheorem{theorem}{Theorem}[section]
\newtheorem{lemma}[theorem]{Lemma}
\newtheorem{proposition}[theorem]{Proposition}

\newtheorem{definition}[theorem]{Definition}
\newtheorem{problem}[theorem]{Problem}

\newtheorem{assumption}[theorem]{Assumption}

\usepackage{eso-pic}
\AddToShipoutPictureBG*{%
\AtPageUpperLeft{%
\setlength\unitlength{1in}%
\hspace*{\dimexpr0.5\paperwidth\relax}
\makebox(0,-1)[c]{
\begin{tabular}{c c}
To be presented at the 64th IEEE Conference on Decision and Control 2025.
Uploaded to ArXiv December 4, 2025. \\
\end{tabular}}}}

\AddToShipoutPictureBG*{%
\AtPageUpperLeft{%
\setlength\unitlength{1in}%
\hspace*{\dimexpr0.5\paperwidth\relax}
\makebox(0,-21.0)[c]{
\footnotesize
\begin{tabular}{c c}
© 2025 IEEE. Personal use of this material is permitted. Permission from
IEEE must be obtained for all other uses, in any \\
current or future media, including reprinting/republishing
this material for advertising or promotional purposes, creating new \\
collective works, for resale or redistribution to servers or lists,
or reuse of any copyrighted component of this work in other works.
\end{tabular}}}}

\def\BibTeX{{\rm B\kern-.05em{\sc i\kern-.025em b}\kern-.08em
    T\kern-.1667em\lower.7ex\hbox{E}\kern-.125emX}}
\begin{document}
\title{\LARGE \bf
Value Function Approximation for Nonlinear MPC: \\ Learning a Terminal Cost Function with a Descent Property
}

\author{T.M.J.T. Baltussen, C.A. Orrico, A. Katriniok, W.P.M.H. Heemels, D. Krishnamoorthy % <-this % stops a space
% \thanks{*This work was not supported by any organization}% <-this % stops a space
\thanks{All authors are with the Control Systems Technology Section, Eindhoven University of Technology, the Netherlands. Dinesh Krishnamoorthy is also with the Dept. of Engineering Cybernetics at the Norwegian University of Science \& Technology, Norway. Email:
{\tt t.m.j.t.baltussen@tue.nl}}}%
% {\tt\{t.m.j.t.baltussen, c.a.orrico, a.katriniok, m.heemels, d.krishnamoorthy\}@tue.nl}}}%

\maketitle

\thispagestyle{empty}

\begin{abstract}
We present a novel method to synthesize a terminal cost function for a nonlinear model predictive controller (MPC) through value function approximation using supervised learning.
Existing methods enforce a descent property on the terminal cost function by construction, thereby restricting the class of terminal cost functions, which in turn can limit the performance and applicability of the MPC. We present a method to approximate the true cost-to-go with a general function approximator that is convex in its parameters, and impose the descent condition on a finite number of states. Through the \textit{scenario approach}, we provide probabilistic guarantees on the descent condition of the terminal cost function over the continuous state space.
We demonstrate and empirically verify our method in a numerical example.
By learning a terminal cost function, the prediction horizon of the MPC can be significantly reduced, resulting in reduced online computational complexity while maintaining good closed-loop performance.
\end{abstract}

\section{Introduction}
\label{sec:Intro}
While model predictive control (MPC) is a powerful method for controlling nonlinear constrained systems, a central challenge in deploying MPC in real-time applications lies in balancing computational efficiency with closed-loop performance.
Firstly, the complexity of the online optimization problem remains a significant challenge for real-time implementation.
Shrinking the prediction horizon can significantly reduce the computational complexity of the MPC.
However, this requires an appropriate terminal cost function that quantifies the \textit{cost-to-go} to obtain acceptable performance. Furthermore, in order to guarantee stability, the closed-loop cost should decrease between consecutive time steps which imposes a \textit{descent condition} on the terminal cost function.
For linear time-invariant systems, this condition is satisfied by using the terminal cost from the infinite-horizon linear quadratic regulator (LQR) problem \cite{rawlings2017model}.
Finding a terminal cost function for nonlinear systems that satisfies this condition is generally a non-trivial task that we aim to address in this work.
Secondly, MPC is increasingly used to tackle complex control tasks by generating optimization-based policies that reflect a local form of Bellman's principle of optimality. However, standard formulations often lack the expressiveness needed to approximate globally optimal behavior, particularly in nonlinear or uncertain settings \cite{banker2025localgloballearninginterpretablecontrol}. To address this, we propose leveraging a richer class of, potentially nonconvex, terminal cost functions to guide the MPC towards solutions that capture more \textit{global} notion of optimality, thereby improving closed-loop performance.

We recognize two streams in the MPC literature, namely explicit and implicit MPC techniques.
In implicit, or \textit{standard} MPC an optimal control problem (OCP) is solved at every time step to compute the control action, which constitutes a control policy. As discussed above, the computational complexity of implicit MPC can be prohibitive.
In some cases, the MPC policy can be computed analytically, which yields an explicit MPC policy \cite{alessio2009survey}.
Alternatively, available data can be used to directly approximate or \textit{learn} the closed-loop MPC policy, e.g. \cite{Storace_PWA_2011, GenLu_CTA12a, Allgower_Policy_2018,Krishnamoorthy_Policy_2021}.
This approximated policy replaces the online optimization, directly mapping a current state to a control input \cite{mesbah2022fusion}.
Despite the potential of machine learning methods,
providing closed-loop guarantees for explicit MPC remains a significant challenge.
In addition, the flexibility of explicit MPC is limited \cite{grancharova2012explicit}, as any changes to the system during its deployment can affect the controller's performance and safety, undermining the advantages of implicit MPC.

\begin{figure}[t]
    \centering
    \vspace{-1mm}\includegraphics[width=0.65\linewidth]{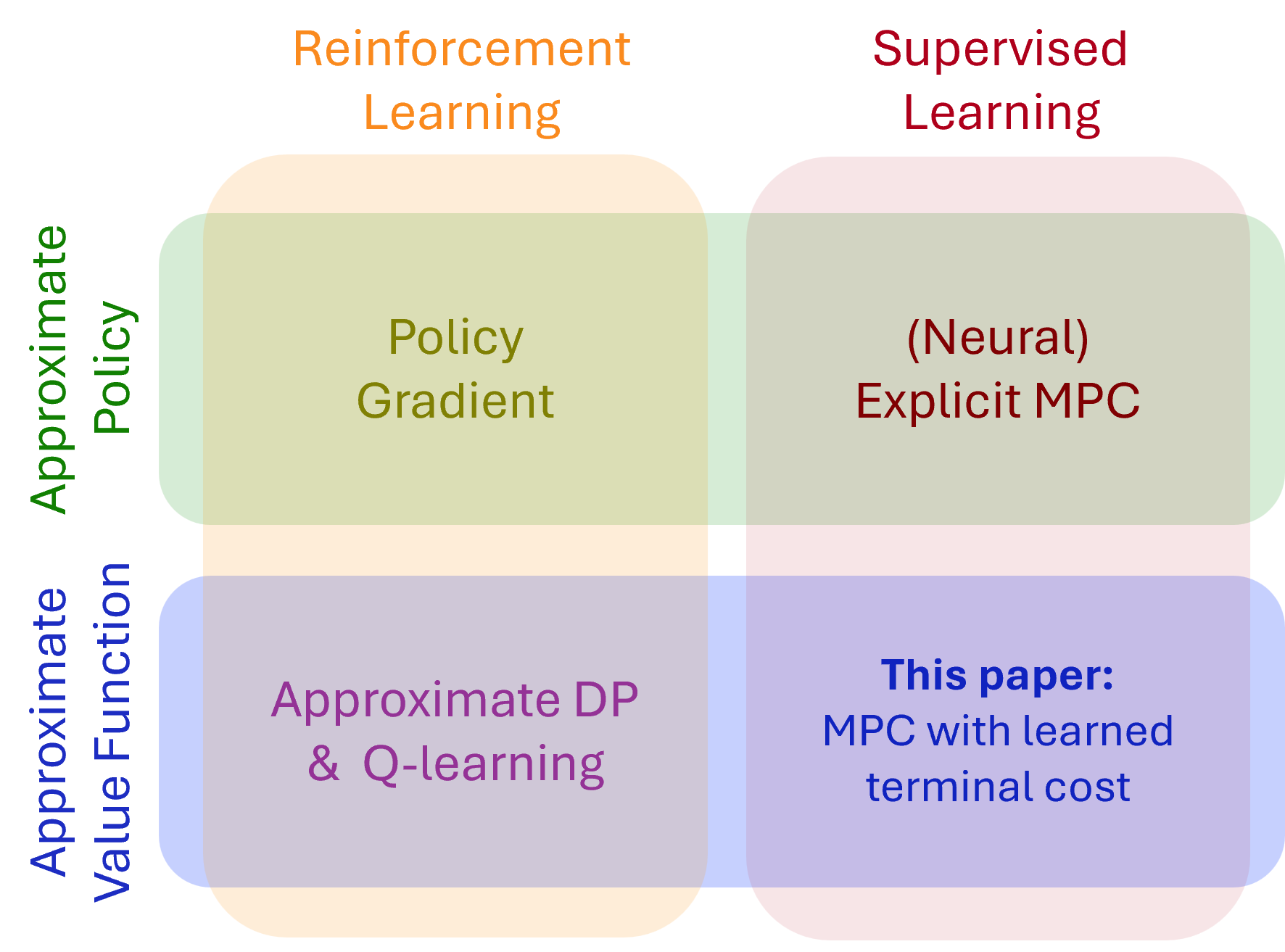}
    \vspace{-1mm}
    \caption{Position of the proposed method against related methods.}
    \label{fig:diagram}
    \vspace{-6mm}
\end{figure}

In this paper, we propose a data-driven method to synthesize a terminal cost function for the OCP of an implicit MPC.
We leverage theoretical results from scenario optimization to give \textit{a priori} probabilistic guarantees on a stabilizing descent property of the learned terminal cost function over the state space.
The proposed method enables a reduction of the horizon length, and hence the online computational complexity of the MPC. In addition, it enables the use of more expressive terminal cost functions that can improve the closed-loop performance when an MPC is used as a locally optimal control policy.
In summary, this paper aims to address these two closely related problems through a data-driven synthesis method for more expressive terminal cost functions that 1) enable a reduction in the MPC horizon length and 2) enable more \textit{global} reasoning when MPC is used as a local control policy, while enforcing and preserving stabilizing properties in a probabilistic manner.
This work lays the foundation for a methodology that aims to reduce the complexity and improve the closed-loop performance of MPC while supporting system-theoretic analysis.

\subsection{Related Work}
A broad class of reinforcement learning (RL) methods addresses challenges closely related to MPC \cite{mesbah2022fusion}. These methods can be categorized in policy approximation and value function approximation methods, cf. Fig. \ref{fig:diagram}.
Learning the terminal value/cost function, or the \textit{cost-to-go} function has been predominantly studied in the context of approximate dynamic programming (ADP) and RL \cite{bertsekas2019reinforcement}, using closed-loop data. That is, the parameters $\theta$ of a parametric function $\mathcal{V}(x; \theta)$ are updated with systematic methods such as Q-learning or temporal difference learning using closed-loop data obtained from different episodes, e.g., using a simulator or system-in-the-loop \cite[Sec. V]{mesbah2022fusion}.
However, the use of supervised learning to learn the cost-to-go function offline for implicit MPC has received little attention and is the focus of this paper.
For a discussion on \textit{value function-augmented} MPC and its relation to RL, we refer the reader to \cite{banker2025localgloballearninginterpretablecontrol}.

The idea of learning a quadratic terminal cost function for linear MPC from optimal state-input data was first presented in \cite{keshavarz2011imputing}.
The problem of learning a quadratic terminal cost function for linear parameter varying (LPV) systems is considered in \cite{abdufattokhov2021learning} and \cite{abdufattokhov2024learning}, where the cost-to-go matrix is a state-dependent matrix whose elements are the outputs of a feedforward neural network.
A neural network-based value function is considered in \cite{mittal2021neurallyapunovmodelpredictive}, however, its descent condition is satisfied by assumption.
In \cite{pauli_learning_2024}, the input-to-state stability with respect to approximation errors in learned value functions is analyzed.
However, simply minimizing the (point-wise) error of the approximated value function does not necessarily lead to the best performance. The value function is intended to steer the MPC. 
In particular, descent of the value function along the closed-loop trajectory is essential for good closed-loop performance.

By learning an appropriate terminal cost function, the MPC horizon can be shortened to reduce the computational complexity of the MPC, with little or no sacrifice in performance \cite{keshavarz2011imputing,abdufattokhov2021learning, abdufattokhov2024learning}.
Despite recent developments, achieving good performance with a short prediction horizon remains an open problem \cite{abdufattokhov2024learning}.
In addition, current methods are limited to (state-dependent) quadratic terminal cost functions, which can introduce conservatism and limit the applicability of this approach.
These methods use terminal cost functions that are convex \textit{by construction}.
By considering a broader class of terminal cost functions, we aim to reduce conservatism, extend the class of systems that we can stabilize using MPC and improve the closed-loop performance of the MPC.

\subsection{Contributions}
We propose a novel supervised learning method to synthesize the terminal cost function of the MPC from a data set of (expert) demonstrations.
This work is motivated by two key challenges: designing terminal cost functions that allow for shorter prediction horizons while maintaining stability, and incorporating richer, potentially nonconvex, terminal cost functions to enhance global decision-making in locally optimal MPC policies.
In this paper, we aim to extend the class of stabilizing terminal cost functions for nonlinear MPC by allowing nonconvex terminal cost functions that are constrained to satisfy the stabilizing descent condition at a finite number of states.
We leverage the scenario approach to provide a probabilistic certificate for the stabilizing condition of the learned terminal cost function over the continuous state space.
To this end, we present the following contributions.
\begin{enumerate}
    \item We present a supervised learning method that enforces the stabilizing descent condition of the terminal cost function on a finite number of points in the state space.
    \item We provide probabilistic guarantees that the learned, nonconvex terminal cost function has a stabilizing descent property everywhere in the state space, except for a region that can be made arbitrarily small.
    \item We demonstrate the proposed method in a numerical example, showing that the horizon length can be significantly reduced with minimal loss in performance compared to an MPC with a long prediction horizon.
    \item As opposed to direct policy learning, we numerically demonstrate the flexibility and robustness of implicit value function-augmented MPC.
\end{enumerate}

\section{Problem Formulation}
\subsection{Stabilization with Nonlinear MPC}
\label{sec:Prob_Stab}
Consider the constrained, discrete-time nonlinear system
\begin{equation} \label{eq:nonlinsystem}
    x_{t+1} = f(x_t, u_t)
\end{equation}
where $x_t \in \mathbb{X} \subset \mathbb{R}^n$ is the system state and $u_t \in \mathbb{U} \subset \mathbb{R}^m$ denotes the control input at time $t \in \mathbb{N} := \{0,1,2, \dots \}$, $\mathbb{X}$ and $\mathbb{U}$ are compact sets that contain the origin, and $f : \mathbb{X} \times \mathbb{U} \rightarrow \mathbb{X}$ is continuously differentiable over $\mathbb{X}$ and $\mathbb{U}$ and $f(0,0) = 0$.
We aim to control system \eqref{eq:nonlinsystem} using standard MPC, wherein we solve the OCP $P_N(x_t)$
\begin{subequations}
\label{eq:MPC}
\begin{align}
\hspace{-1mm}
V(x_t) :&= 
\min _{\mathbf{u}} \, \sum_{k=0}^{N-1} \ell\left(x_{k \mid t}, u_{k \mid t}\right) + \mathcal{V}(x_{N \mid t}; \theta), \\
\text { s.t. } x_{k+1 \mid t} & =f(x_{k \mid t}, u_{k \mid t} ), \quad \hspace{-2mm} \forall k \in \{0,1, \dots, N - 1\}, \\
u_{k \mid t} &\in \mathbb{U}, \quad \hspace{13.5mm} \forall k \in \{0,1, \dots, N - 1\}, \\
x_{k \mid t} & \in \mathbb{X}, \quad \hspace{13.5mm} \forall k \in \{1,2, \dots, N\}, \\
x_{0 \mid t} & = x_t,
\end{align}
\end{subequations}
with a terminal cost function $\mathcal{V}(\cdot;\theta) : \mathbb{R}^n \rightarrow \mathbb{R}_{\geq 0}$ parameterized by $
\theta \in \mathbb{R}^d$.
The stage cost function $\ell : \mathbb{X} \times \mathbb{U} \rightarrow \mathbb{R}_{\geq 0}$, where $\ell(0,0) = 0$, is typically taken to be a quadratic function of the state and input.
The value function $V$ denotes the minimum of the MPC cost function:
\begin{equation}
    J(\mathbf{x}_t, \mathbf{u}_t) := \sum_{k=0}^{N-1} \ell\left(x_{k \mid t}, u_{k \mid t}\right) + \mathcal{V}(x_{N \mid t}; \theta),
\end{equation}
where $\mathbf{x}_t = (x_{0 \mid t}, x_{1 \mid t}, \ldots, x_{N \mid t})$ is a state sequence and $\mathbf{u}_t = (u_{0 \mid t}, u_{1 \mid t}, \ldots, u_{N-1 \mid t})$ an input sequence computed at time $t \in \mathbb{N}$.
The resulting MPC control law is given by $u_t = \kappa_N(x_t) := u^*_{0 \mid t}$, with $\mathbf{u}_t^* = (u^*_{0 \mid t}, u^*_{1 \mid t}, \ldots, u^*_{N-1 \mid t})$ being the minimizer of \eqref{eq:MPC}, assuming that this exists.
Let $\mathcal{X} \subseteq \mathbb{X}$ denote the feasible set of states of MPC \eqref{eq:MPC}. We assume that the MPC is feasible under the control law $\kappa_N$, i.e. we assume that the MPC is recursively feasible.

\begin{assumption}
    When $P_N(x_t)$ is feasible for a state $x_t$ at time step $t$, i.e. $x_t \in \mathcal{X}$, then $P_N(x_{t+1})$ is feasible at time $t+1$ for the state $x_{t+1} = f(x_t, \kappa_N(x_t))$ such that $x_{t+1} \in \mathcal{X}$.
\end{assumption}
\vspace{0.5em}

\noindent
Recursive feasibility is an essential property of an MPC.
In many practical cases, soft-constraints can be used to ensure recursive feasibility of the MPC \cite{rawlings2017model}. While addressing recursive feasibility remains highly relevant, in this paper we focus on the stability properties of the MPC.
The closed-loop stability of MPC is typically analyzed using the value function $V$.
In order to provide asymptotic stability of the origin under the closed-loop MPC, typically three main conditions on $V$ are required, namely an upper bounding comparison function and a lower bounding comparison function on $V$, and a \textit{descent condition} on $V$ along the closed-loop trajectories of the MPC, see, e.g., \cite{rawlings2017model}.\\

\begin{assumption}
\label{ass:lowerbound}
    There exist $\alpha_1, \alpha_2 \in \mathcal{K}_\infty$\footnote{A continuous function $\alpha : [0, \infty) \rightarrow [0, \infty)$ is said to be of class $\mathcal{K}_\infty$ if it is strictly increasing, $\alpha(0) = 0$ and $\lim\limits_{s \rightarrow \infty} \alpha(s) = \infty$.} such that
    \begin{equation}
        \alpha_1(\|x\|) \leq V(x) \leq \alpha_2(\|x\|) \quad \forall x \in \mathcal{X}.
    \end{equation}
\end{assumption}
\vspace{0.5em}

The upper and lower bounding comparison functions $\alpha_1, \alpha_2 \in \mathcal{K}_\infty$ of $V$ are typically obtained from the quadratic stage cost $\ell$ and the compactness of $\mathbb{X}$ and $\mathbb{U}$.
The MPC value function $V$ serves as a valid Lyapunov function, if the terminal cost function $\mathcal{V}$ satisfies a basic stability condition \cite{rawlings2017model} that we define as follows.\\

\begin{definition}\textit{(Basic stability condition)}
    \label{def:stability_cond}
    A terminal cost function $\mathcal{V}$ is said to satisfy the basic stability condition if for all $x \in \mathcal{X}$ there exists an admissible $u \in \mathbb{U}$ such that:
    \begin{align}
    \label{eq:descentcond}
        \mathcal{V}(f(x,u)) & - \mathcal{V}(x) \leq - \ell(x,u).
    \end{align}
\end{definition}
\vspace{0.5em}

\noindent
This condition \eqref{eq:descentcond} is typically assumed to hold by design within a control invariant terminal set $\mathbb{X}_T \subseteq \mathcal{X}$. The terminal state $x_{N \mid t}$ is then constrained to lie in $\mathbb{X}_T$.
For nonlinear systems, a quadratic terminal cost can be used to satisfy the descent condition \eqref{eq:descentcond} for a local linearization around the origin \cite{rawlings2017model}.
However, a terminal constraint set $\mathbb{X}_T$ can restrict the size of $\mathcal{X}$, in particular when the prediction horizon is very short.
In order to increase the domain of attraction, we seek to satisfy \eqref{eq:descentcond} over the feasible set of states $\mathcal{X}$.
Rather than restricting $\mathcal{V}(x;\theta)$ to a specific class of functions, e.g. quadratic in $x$, that ensures stability on a reachable subset $\mathbb{X}_T \subseteq \mathcal{X}$,
we consider $\mathcal{V}(\cdot;\theta) : \mathbb{R}^n \rightarrow \mathbb{R}_{\geq 0}$ to be a general function approximator that is linear in its parameters.
For example, $\mathcal{V}(x;\theta)$ can be a linear combination of nonlinear basis functions $\phi(x)$:
\begin{equation}
\label{eq:cost_fun}
    \mathcal{V}(x;\theta) = \theta^\top \phi(x).
\end{equation}
This allows us to capture nonlinear, nonconvex functions using function approximators that are convex in their parameters.
In many practical scenarios, data samples of the state, input and value function are available that can be leveraged to synthesize a terminal cost function. These data can stem from a baseline controller or (human) expert demonstrations that can be computed offline.
We can enforce the \textit{basic stability condition} from Definition \ref{def:stability_cond} by imposing it as an explicit constraint in the synthesis problem.
To this end, we address the following problem.\\

\begin{problem}
    \label{def:Problem_Def}
    Find a parameter $\theta \in \mathbb{R}^d$ such that the parametric terminal cost function $\mathcal{V}(x;\theta)$ of the MPC \eqref{eq:MPC} satisfies the basic stability condition from Definition \ref{def:stability_cond} for almost all states with high probability.
\end{problem}

\section{Scenario-based Synthesis}
In this section, we first define a convex synthesis problem to address Problem \ref{def:Problem_Def}. Secondly, we discuss the scenario approach that we will leverage in subsequent sections. We then propose a randomized synthesis program to efficiently solve the synthesis problem based on uniform samples from the state space. Lastly, we discuss the sampling method that we use for the randomized synthesis problem.

\subsection{Value Function Learning}
We define the following convex synthesis problem $S$, enforcing a descent condition of $\mathcal{V}(x; \theta)$ over the set $\mathcal{X}$: 
\begin{subequations}
\label{eq:Synthesis_Prob}
\begin{align}
    S \, &: \, 
    \min _{\theta \in \mathbb{R}^d} \mathcal{G}\bigl(\theta\bigr) \\
    \text{s.t. } &\mathcal{V}\bigl(f\bigl(x, \kappa(x)\bigr) ; \theta\bigr) - \mathcal{V}\bigl(x ; \theta\bigr) \leq -\ell \bigl(x, \kappa(x) \bigr), 
    \label{eq:Descent_Constraint} \notag \\
    & \forall  x \in \mathcal{X}.
    \end{align}
\end{subequations}
\vspace{-5mm}

\noindent
where $\kappa : \mathcal{X} \rightarrow \mathbb{U}$ denotes an arbitrary, but admissible control policy.
Here, we try to find a parameter $\theta \in \mathbb{R}^d$ that ensures that Definition \ref{def:stability_cond} is satisfied for some admissible input $u = \kappa(x) \in \mathbb{U}$ for all states $x \in \mathcal{X}$, while minimizing some metric $\mathcal{G}$ that is convex in $\theta$.
The synthesis problem $S$ has a finite number of optimization variables, but an infinite number of constraints since $\mathcal{X}$ is a continuous space. Such an optimization problem is referred to as a semi-infinite program \cite{Campi_CVX_2008} and is generally intractable.
Instead, we can extract samples from the continuous set $\mathcal{X}$ and solve an approximation of $S$. 
The question that then arises is what is the probability that $\eqref{eq:Descent_Constraint}$ holds for any state $x \in \mathcal{X}$, given \eqref{eq:Descent_Constraint} holds for a finite number of states?
To answer this question, we leverage results from the 
\textit{scenario approach}.

\subsection{Preliminaries on Scenario Optimization}
Robust optimization problems are used for decision making under uncertainty, and account for all possible uncertainties in some set $\Delta$.
The scenario approach \cite{Campi_CVX_2008} is a general randomized decision-making approach that uses a finite number of samples from this set $\Delta$ such that the resulting problem can be solved at a relatively low computational cost \cite{Campi_CVX_2008}.
It has been successfully applied in various control problems such as stochastic MPC, interval predictor models and linear matrix inequalities \cite{campi2009scenario}.
We refer the reader \cite{campi2009scenario} for a comprehensive introduction.
Consider the following general \textit{scenario program}:
\begin{subequations}
\label{eq:Random_Prog}
\begin{align}
    RP_M \, &: \, 
    \min _{\theta \in \mathbb{R}^d} c^\top \theta \\
    \text{s.t. } & \theta \in \bigcap\limits_{i \in \{1, \dots, M \} }  \Theta_{\delta^{(i)}},
    \end{align}
\end{subequations}
\vspace{-2mm}

\noindent
where $M$ samples (or scenarios) $\delta^{(1)}, \delta^{(2)}, \dots, \delta^{(M)}$ are randomly extracted from $\Delta^M$ according to a distribution $\mathbb{P}_\delta$.
Each sample $\delta$ defines a constraint $\Theta_\delta$. The collection of $M$ constraints is then simultaneously enforced in $RP_M$.
Note that the dependence of $\Theta_\delta$ on $\delta$ can be highly nonlinear.
In order to assess the constraint violation in the complete set $\Delta$ by the solution $\theta^*_M$ of the scenario program $RP_M$, the following notion of \textit{violation probability} is central \cite{Campi_CVX_2008}.\\

\vspace{-1mm}
\begin{definition}
\label{def:viol}
    The \emph{violation probability} of a given $\theta \in \mathrm{\Theta}$ is defined as $\mathcal{P}(\theta)=\mathbb{P}_\delta\left\{\delta \in \Delta: \theta \notin \mathrm{\Theta}_\delta\right\}$.
\end{definition}
\vspace{0.5em}

The violation probability $\mathcal{P}(\theta)$ is the probability of drawing a new sample $\delta \in \Delta$ for which the solution $\theta$ does not satisfy the constraint $\Theta_\delta$ associated with the realization $\delta$. 
Note that the violation probability $\mathcal{P}(\theta^*_M)$ is a random variable, since it depends on the random extractions $\delta^{(1)}, \delta^{(2)}, \dots, \delta^{(M)}$. The scenario approach provides a confidence bound $\beta$ on the violation probability $\epsilon$ of $\theta^*_M$.\\

\vspace{-1mm}
\begin{lemma}\textit{($\epsilon$-$\beta$ Result \cite{Campi_CVX_2008})} \\
\label{lem:epsBetaResult}%
    Under the existence and uniqueness of the solution to $RP_M$:
    \begin{equation}
        \mathbb{P}^M\left\{\mathcal{P}\left(\theta_M^*\right)>\epsilon\right\} \leq \sum_{i=0}^{d-1}\binom{M}{i} \epsilon^i(1-\epsilon)^{M-i} = \beta,
    \end{equation}
    where $\epsilon$ denotes the violation parameter and $\beta$ denotes the confidence parameter.
\end{lemma}
\vspace{0.5em}
\noindent
Here, the binomial coefficient of $M$ and $i$ is depicted by $\binom{M}{i}$.
The $\epsilon$-$\beta$ result (Lemma \ref{lem:epsBetaResult}) shows that the violation probability of the optimal solution to $RP_M$ $(\theta^*_M)$ is bounded by a binomial distribution that depends on the number of samples $M$ and the dimension of $\theta$.
Furthermore, the result tells us how many samples $M$ are required to bound the violation probability by $\epsilon$ with a desired level of confidence $1-\beta$. For details on the scenario approach, see \cite{Campi_CVX_2008, campi2009scenario}.

\subsection{The Randomized Synthesis Program}
\label{sec:Random_Synthesis}
Subsequently, we will leverage the results from the scenario approach to approximate the semi-infinite program $S$ \eqref{eq:Synthesis_Prob} by a randomized synthesis program $S_M$ based on $M$ samples from $\Delta$.
Suppose we are given an independent and identically distributed (i.i.d.) data set $\mathcal{D} = \left\{\delta^{(i)} \right\}_{i = 1}^M$ of state-input-value samples from the sample space $\Delta$:
\begin{equation}
    \delta^{(i)} = \left(x^{(i)}, u^{(i)}, \mathcal{J}\bigl(x^{(i)} \bigr) \right)
\end{equation}
with $x^{(i)} \in \mathcal{X}$, $u^{(i)} \in \mathbb{U}$. The input samples are assumed to be generated from a baseline controller $\kappa : \mathcal{X} \rightarrow \mathbb{U}$:
\begin{equation}
    \label{eq:Demo_policy}
    u^{(i)} = \kappa\bigl(x^{(i)}\bigr).
\end{equation}
Here, the value samples $\mathcal{J}\bigl(x^{(i)}\bigr)$ originate from a (possibly unknown) value function $\mathcal{J} : \mathbb{R}^{n_x} \rightarrow \mathbb{R}_{\geq0}$ that quantifies the cost-to-go.
For example, when approximating a long horizon MPC, $\kappa$ denotes the long horizon MPC policy, and $\mathcal{J}$ denotes the value function of this long horizon MPC.
Through the scenario approach, we do not require access to these functions, but just to samples of these functions.
Note that if we do not have data from $\mathcal{J}$, the proposed method can still be applied by learning for state-input data through inverse optimization, as is done, for example, in \cite{keshavarz2011imputing}.

We then construct the following scenario program $S_M$ with the \textit{basic stability condition} from Definition \ref{def:stability_cond} enforced on $\mathcal{V}(x; \theta)$ at the $M$ sampled states to find $\theta^*_M$.
\begin{subequations}
\begin{align}
    S_M \, &: \, 
    \min _{\theta \in \mathbb{R}^d} \sum_{i} \mathcal{G}\bigl(\mathcal{V}\bigl(x^{(i)} ; \theta\bigr), \mathcal{J}\bigl(x^{(i)}\bigr),\theta\bigr) \\
    \text{s.t. } \mathcal{V} &\bigl(f\bigl(x^{(i)},u^{(i)}\bigr) ; \theta\bigr) - \mathcal{V}\bigl(x^{(i)} ; \theta\bigr) \leq -\ell \bigl(x^{(i)}, u^{(i)} \bigr) \label{eq:Scenario_Constraint} \notag \\
    & \forall  i \in \{1, \dots, M\},
    \end{align}
    \label{eq:Scenario_Program}
\end{subequations}
where the constraint \eqref{eq:Scenario_Constraint} defines $\Theta_{\delta^{(i)}}$.
Let us consider the synthesis problem in \eqref{eq:Synthesis_Prob}, and assume that $\Delta = \mathcal{X} \times \mathbb{U} \times \mathbb{R}_{\geq0}$ is equipped with a probability distribution $\mathbb{P}_\delta$ that is uniform over $\mathcal{X}$. Moreover, we assume that the class of $\mathcal{V}(x;\theta)$ is sufficiently rich such that the solution to the scenario program $\theta^*_M$ exists and is unique. By synthesizing a terminal cost function through $S_M$ \eqref{eq:Scenario_Program}, we obtain a probabilistic certificate for the violation probability of $\theta^*_M$ to the original synthesis problem $S$ \eqref{eq:Synthesis_Prob}.

\subsection{Uniform Sampling}
For the terminal cost synthesis problem \eqref{eq:Scenario_Program} we employ a uniform sampling method. 
To this end, we assume that states $x^{(1)}, x^{(2)}, \dots, x^{(M)}$ from the scenarios $\delta^{(1)}, \delta^{(2)}, \dots, \delta^{(M)}$ are uniformly sampled from $\mathcal{X}^M$, where $\mathcal{X}^M$ denotes the product space of sampling domain $\mathcal{X}$.
The input and cost samples can take an arbitrary distribution, as this does not affect our analysis over the set $\mathcal{X}$.
The motivation for this sampling approach is detailed below.

\subsubsection{Volume of the violation set}
Firstly, let us denote the violation set of the scenario program by $\Delta_\epsilon(\theta^*_M) := \{\delta \in \Delta : \theta \notin \Theta_\delta\}$. By taking uniform samples from $\mathcal{X}$, we obtain a bound on the volume of the \textit{violation set} denoted by $\Delta_\epsilon(\theta) \subset \Delta$ \cite[Ch. 6]{tempo2013randomized} on the feasible set of states $\mathcal{X}$. Under uniform sampling, the violation probability determines the ratio between the volume of the violation set and the volume of the domain.
This certificate directly follows from the definition of the violation probability and Lemma \ref{lem:epsBetaResult}. \vspace{0.5em}
\begin{proposition}
    \label{prop:fraction}
    Suppose that $\delta \in \Delta$ adheres to a uniform distribution $\mathbb{P}_\delta$. Then, the ratio of the violation set $\Delta_\epsilon(\theta^*_M)$ and the sample space $\Delta$ is bounded by $\epsilon$.
\end{proposition}
\begin{proof}
    From the $\epsilon$-$\beta$ result (Lemma \ref{lem:epsBetaResult}) we have
    \begin{equation}
        \mathcal{P}(\theta^*_M) \leq \epsilon,
    \end{equation}
    with probability $1-\beta$.
    Let $p(\delta)$ denote the uniform probability density function of $\delta$. Then,
\begin{align}
    &\mathcal{P}(\theta^*_M) = \int_{\Delta_\epsilon}p(\delta ) \, \mathrm{d} \delta = \int_{\Delta_\epsilon}\frac{1}{\text{Vol}(\Delta)} \, \mathrm{d} \delta  \leq \epsilon \\
    &\implies \text{Vol}(\Delta_\epsilon) = \int_{\Delta_\epsilon} \mathrm{d} \delta \leq \epsilon \text{Vol}(\Delta) \\
    &\implies \frac{\text{Vol}(\Delta_\epsilon)}{\text{Vol}(\Delta)} \leq \epsilon,
\end{align}
    with probability no smaller than $1-\beta$.
\end{proof}
\vspace{0.5em}

\noindent
Proposition \ref{prop:fraction} implies that, under uniform sampling of $\mathcal{X}$, the descent condition of the value function from $S_M$ is satisfied everywhere, but at most an $\epsilon$-volume of the feasible set of states $\mathcal{X}$, i.e., the violation set.
Hence, by uniform sampling of $\mathcal{X}$ we can control the size of the violation set and we can have high confidence that the learned function $\mathcal{V}$ is a suitable terminal cost function for the MPC.

\section{Main Result}
\subsection{Probabilistic Certificate of the Descent Condition}
The $\epsilon$-$\beta$ result from the scenario approach \cite{Campi_CVX_2008} certifies that, under a uniform distribution, the basic stability condition of the terminal cost function is satisfied for all states except for an $\epsilon$-fraction with high probability.\\

\begin{theorem}
\label{th:Stab}
    \textit{(Main Result)} \newline
    The learned terminal cost function $\mathcal{V}(x, \theta^*_M)$ satisfies the basic stability condition from Definition \ref{def:stability_cond} for all $x \in \mathcal{X}$ except for at most an $\epsilon$-fraction, with probability no smaller than $1-\beta$.
\end{theorem}
\begin{proof}
Due to Lemma \ref{lem:epsBetaResult}, we obtain a bound on the violation probability of the descent property of the learned terminal cost function $\mathcal{V}$:
\begin{equation}
\label{eq:prob_beta}
    \mathbb{P}\{ \mathcal{V}(f(x, \kappa(x)) ; \theta_M^*) - \mathcal{V}(x ; \theta_M^*) \leq - \ell(x, \kappa(x)) \} >\epsilon,
\end{equation}
with at most probability $\beta$. This probability holds for all $x \in \mathcal{X}$, due to the uniform distribution over $\mathcal{X}$.
By taking the complement of \eqref{eq:prob_beta}, we obtain:
\begin{equation}
\label{eq:probDescentAnyU}
    \mathbb{P}\{ \mathcal{V}(f(x, \kappa(x)) ; \theta_M^*) - \mathcal{V}(x ; \theta_M^*) \leq - \ell(x, \kappa(x)) \} \leq \epsilon,
\end{equation}
with at least probability $1 - \beta$ for all $ x \in \mathcal{X}$.\\

\noindent
Due to Proposition \ref{prop:fraction}, this implies that the violation set is an $\epsilon$-fraction of $\mathcal{X}$ under uniform distribution over $\mathcal{X}$, such that
\begin{equation}
    \mathcal{V}(f(x, \kappa(x)) ; \theta_M^*) - \mathcal{V}(x ; \theta_M^*) \leq - \ell(x, \kappa(x))
\end{equation}
for all $x \in \mathcal{X}$ but at most an $\epsilon$-fraction with probability $1-\beta$. \\

\noindent
Note that this certificate holds for the demonstrating policy $\kappa(x)$ \eqref{eq:Demo_policy}, which is suboptimal to the MPC problem $P_N$ with the learned terminal cost $\mathcal{V}(x, \theta^*_M)$. 
By optimality, there exists an admissible (sub)optimal control input $u\in \mathbb{U}$ such that,
\begin{equation}
\label{eq:probDescentUstar}
    \mathcal{V}(f(x, u) ; \theta_M^*) - \mathcal{V}(x ; \theta_M^*) \leq - \ell(x, u),
\end{equation}
for all $x \in \mathcal{X}$ but at most an $\epsilon$-fraction with probability $1-\beta$.
\end{proof}
\vspace{0.5em}

\noindent
Note that the scenario approach does not provide any guarantees if the $M$ samples are drawn from the $\beta$ volume in $\Delta^M$, accounting for the probability that the samples are drawn from a `bad' set. However, the size of $\beta$ can be made extremely small considering the cheap computational complexity of $\beta$ \cite{Campi_CVX_2008}.
Hence, we conclude that the basic stability condition of the learned terminal cost function $\mathcal{V}(\cdot;\theta^*_M)$ can be enforced everywhere in $\mathcal{X}$, except for a region $\Delta_\epsilon$ that can be made arbitrarily small.

\subsection{Extension to Myopic MPC}
Value function approximation is of great interest for MPC with an extremely short horizon, as short as $N=1$. The performance of such a \textit{myopic} MPC strongly relies on the accuracy of $\mathcal{V}$. 
This problem is closely related to approximate dynamic programming:
\begin{subequations}
\begin{align}
\min_{u_t} \, & \ell(x_t, u_t)+ \mathcal{V}\left(x_{t+1}\right), \\
\text {s.t. } x_{t+1} & = f(x_t, u_t), \\
x_{t+1} & \in \mathbb{X},\\
u_t & \in \mathbb{U},
\end{align}    
\end{subequations}
where $\mathcal{V}$ approximates true value function $\mathcal{J}$.
This closely resembles approximate dynamic programming \cite{bertsekas2019reinforcement}, with the key difference that our approximate cost-to-go function $\mathcal{V}$ is learned offline using supervised learning, without the need for a simulator or a system-in-the-loop, cf. Fig. \ref{fig:diagram}.

\section{Numerical Illustration}
\subsection{Continuous Stirred Tank Reactor}
We now apply the proposed approach on a benchmark continuous stirred tank reactor (CSTR) problem that was also used in the context of MPC policy approximation in \cite{Allgower_Policy_2018, Krishnamoorthy_Policy_2021} and value function learning in \cite{abdufattokhov2024learning}. This problem consists of two states, the scaled concentration and the reactor temperature of a CSTR, denoted by $x_1$ and $x_2$, respectively. The process is controlled using the coolant flow rate $u$. The continuous-time dynamics are given by:
\begin{align}
& \dot{x}_1=(1 / \tau)\left(1-x_1\right)-k x_1 e^{-b / x_2}, \\
& \dot{x}_2=(1 / \tau)\left(x_f-x_2\right)+k x_1 e^{-b / x_2} - a u\left(x_2-x_c\right),
\end{align}
where the model parameters are $\tau=20$, $k=300$, $b=5$, $x_f=$ $0.3947$, $x_c=0.3816$, and $a=0.117$. Furthermore, we have $\mathbb{X}= [0.0632,0.4632] \times[0.4519, 0.8519]$ and $\mathbb{U}=[0,2]$. 
The system is discretized using Euler discretization with a discretization step of $h = 0.1$ [s].
The set point is given by $x^{s p}=[0.2632,0.6519]^{\top}$ with the associated steady state control input $u^{sp} = 0.7853$. For the MPC we apply a change of variables and control the shifted system to the origin. Hence, the stage cost is given by
\begin{equation}    
\ell(x, u)= \left\|x-x^{s p}\right\|^2+10^{-4}\|u - u^{sp}\|^2.
\end{equation}

\begin{figure*}[b!]
    \centering
    \begin{subfigure}{0.26\textwidth} % Adjust width to fit three in a row
        \centering
        \includegraphics[width=\linewidth]{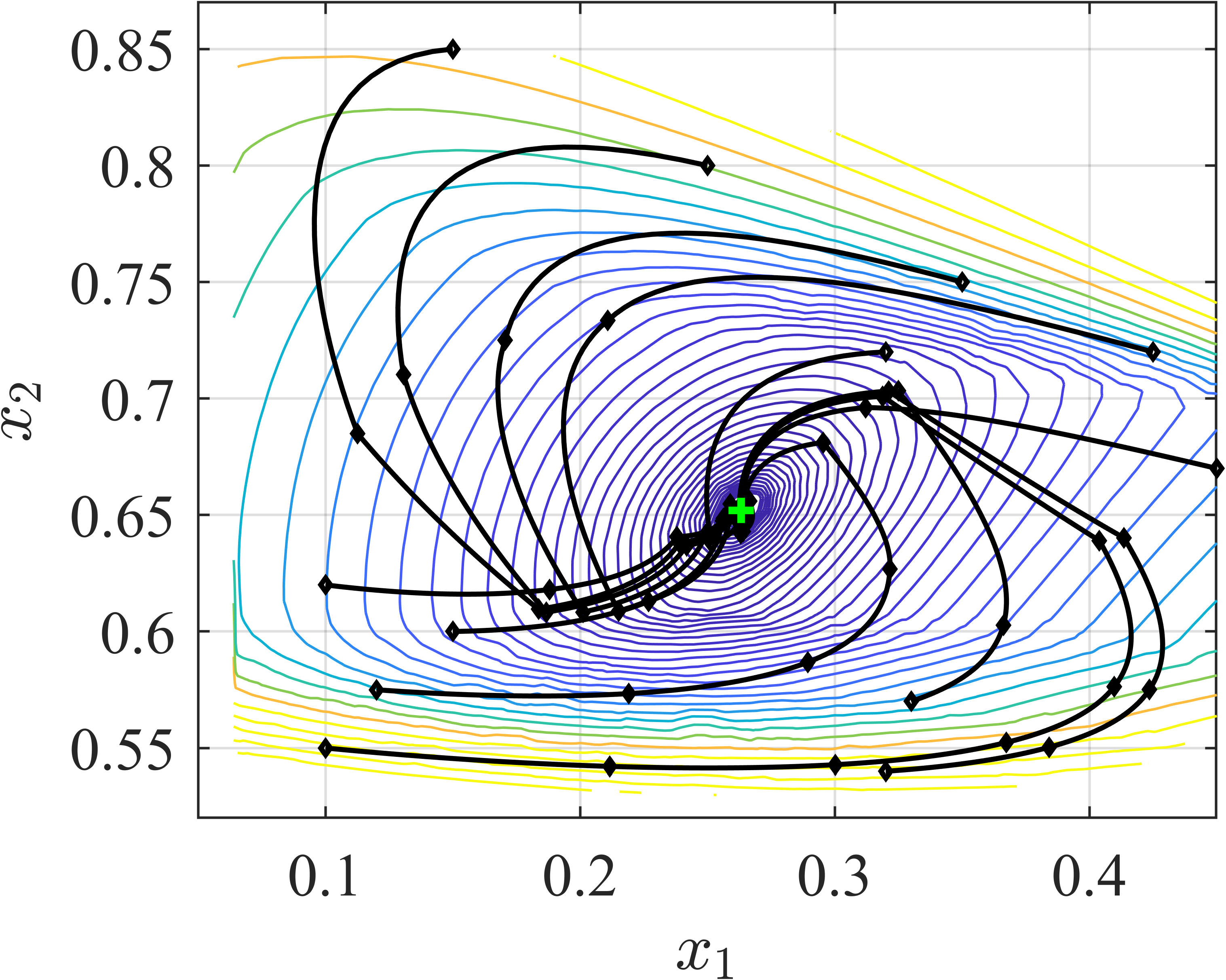}
        \caption{Expert MPC. \\\hspace{\textwidth} }
        \label{fig:cstr_full_nmpc}
    \end{subfigure}
    \hfill
    \begin{subfigure}{0.26\textwidth}
        \centering
        \includegraphics[width=\linewidth]{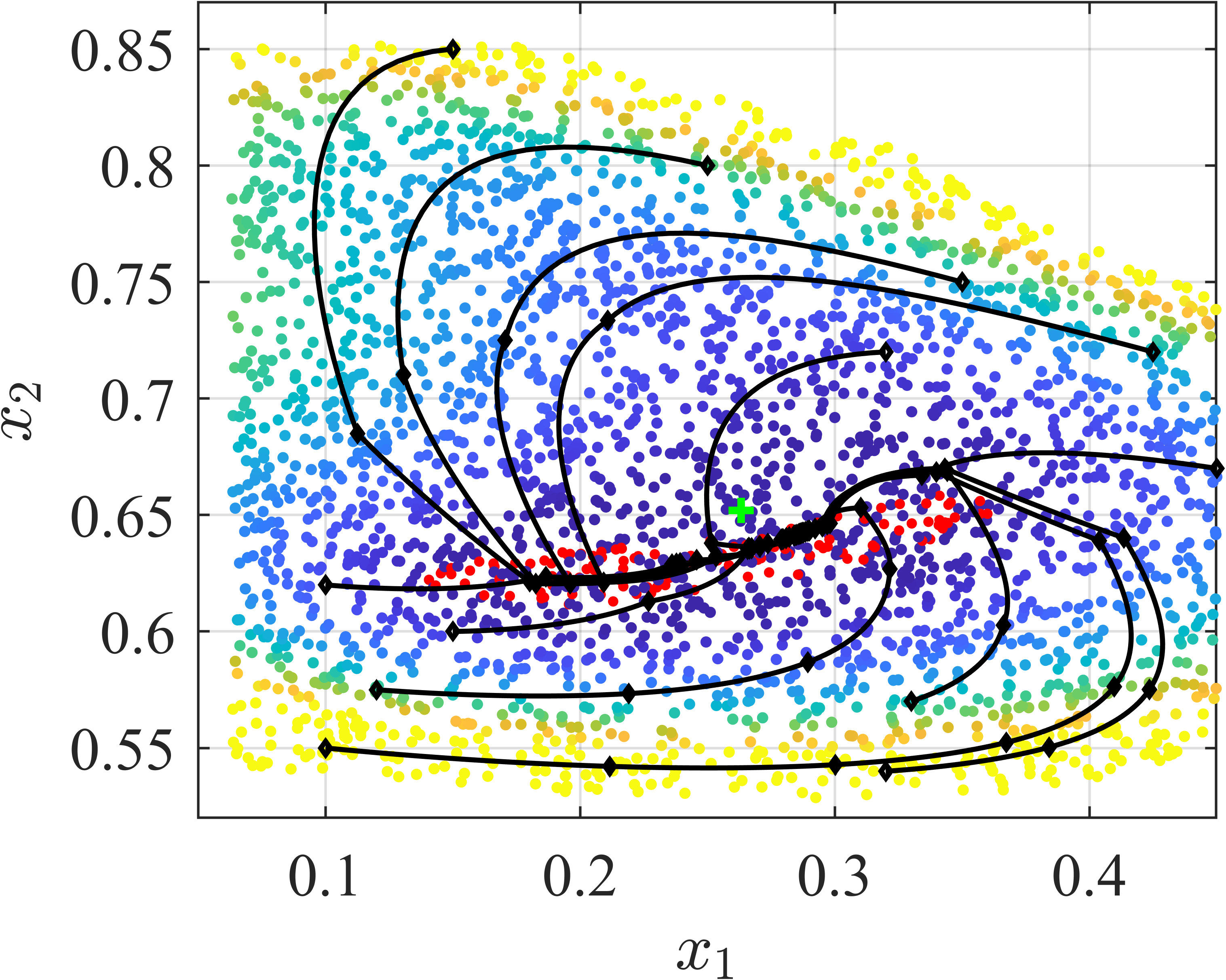}
        \caption{Proposed MPC \textit{without} the \newline descent constraint.}
    \label{fig:cstr_nodesc_viol}
    \end{subfigure}
    \hfill
    \begin{subfigure}{0.44\textwidth}
        \centering
        \includegraphics[width=\linewidth]{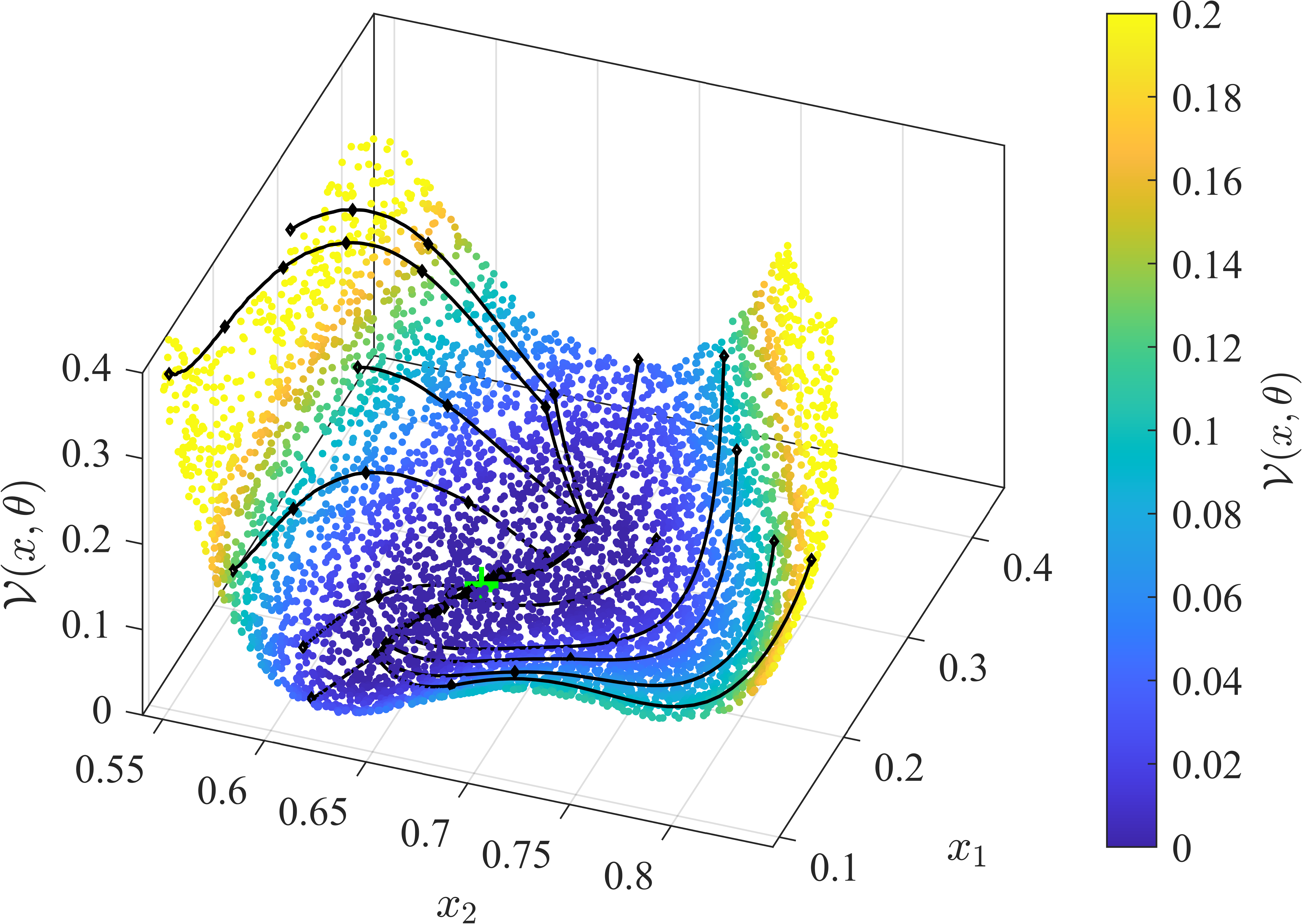}
        \caption{Proposed MPC \textit{with} the descent constraint. \\\hspace{\textwidth}}
    \label{fig:cstr_3d}
    \end{subfigure}
    \caption{Phase plots of twelve closed-loop trajectories of (a) the expert $(T=50)$ MPC and proposed MPC $(N=1)$ from $S_M$ with $(M=2842)$ \text{(b) \textit{without}} and (c) \textit{with} the descent constraint \eqref{eq:Scenario_Descent}. The color-graded (a) level sets depict the true cost-to-go $\mathcal{J}(x)$ and (b-c) the dots depict the learned cost-to-go $\mathcal{V}(x, \theta^*_M)$ evaluated at the training points. Red dots in (b) indicate the violation of the descent condition in the test data set. The green cross indicates the set point of the MPC. The nonconvexity of $\mathcal{V}$ can be observed in (c).}
    \label{fig:nmpc_compare}
\end{figure*}

\subsection{Expert Demonstration}
The synthesis problem relies on a set of expert demonstrations that are generated by a baseline policy $\kappa$. In this example, these training data are generated by an \textit{expert} MPC,
\begin{subequations}
\label{eq:Expert_MPC}
    \begin{align}
        \hspace{-10mm}
        \mathcal{J}(x_t) : = \min_{\mathbf{u}} & \, \sum_{k=0}^{T-1} \left\|x_{k \mid t}-x^{s p}\right\|^2 + 10^{-4}\|u_{k \mid t} - u^{sp}\|^2 \notag \\ & \quad + \left\|x_{T \mid t}-x^{s p}\right\|^2, \\
        \text { s.t. } x_{k+1 \mid t} = & \, f(x_{k \mid t}, u_{k \mid t} ), \hspace{0.5mm} \forall k \in \{0,1, \dots, T - 1\}, \\
        u_{k \mid t} \in & \, \mathbb{U}, \hspace{16mm} \forall k \in \{0,1, \dots, T - 1\}, \\
        x_{k \mid t} \in & \, \mathbb{X}, \hspace{16mm} \forall k \in \{1,2, \dots, T\}, \\
        x_{0 \mid t} = & \, x_t,
        \end{align}
\end{subequations}
designed to stabilize the CSTR to the desired unstable setpoint $(x^{sp}, u^{sp})$ and features a sufficiently long prediction horizon $T >> N$ such that the system is stable without terminal ingredients
\cite[Theorem 3.6]{GRUNE21}, such as a control invariant terminal set and an infinite-horizon LQR cost.
Here, we take $T=50$ and solve the MPC at a sampling time of $T_{s}=3 \mathrm{~[s]}$.
We generate a total of $M_{\mathrm{data}}=9068$ state-input-value tuples $\bigl(x_t, \kappa(x_t), \mathcal{J}(x_t)\bigr)$, by solving \eqref{eq:Expert_MPC} at uniform samples from the feasible set of states $\mathcal{X}$.
In addition, a set $\mathcal{D}_{\text{test}}$ of $M_{\text{test}} = 6226$ state-input-value tuples $\bigl(x_t, \kappa(x_t), \mathcal{J}(x_t)\bigr)$ is generated for validation.

\subsection{Value Function Learning}
We use the data set $\mathcal{D}_{\text{data}} = \{x^{(i)}, u^{(i)}, \mathcal{J}(x^{(i)})\}_{i = 1}^{M_{\text{data}}}$ from the \textit{expert} MPC \eqref{eq:Expert_MPC} as our training data set and solve the following scenario program:
\begin{subequations}
\label{eq:Scenario_Program_CSTR}
\begin{align}
    S_M \, &: \, 
    \min _{\theta \in \mathbb{R}^d} \frac{1}{M} \sum_{i=1}^{M}\left\|\mathcal{V}\left(x_i ; \theta\right)-\mathcal{J}\left(x_i\right)\right\|^2 \\
    \text{s.t. } &\mathcal{V}\bigl(f\bigl(x^{(i)},u^{(i)}\bigr) ; \theta\bigr) - \mathcal{V}\bigl(x^{(i)} ; \theta\bigr) \leq -\ell \bigl(x^{(i)}, u^{(i)} \bigr), \notag \\
    & \forall  i \in \{1, \dots, M\}
    \label{eq:Scenario_Descent},
    \end{align}
\end{subequations}
where $\mathcal{V}\left(x ; \theta\right) = \theta^\top \phi(x)$ \eqref{eq:cost_fun} is a radial basis function network with $d = 75$ neurons, with $\phi : \mathbb{R}^n \rightarrow \mathbb{R}^d$.
We minimize the mean squared error of the terminal cost function while enforcing the descent condition at the sampled states.
Note that the theoretical results only rely on the convexity of the scenario program. However, as mentioned in Section \ref{sec:Random_Synthesis}, we require that the class of terminal cost functions is sufficiently rich such that we have a feasible solution to the synthesis problem, which may be problem-specific.
Furthermore, the performance of the proposed MPC will depend on the tuning of the expert data as the learned terminal cost function will try to `mimic' the expert control policy.
We consider a confidence of $\beta = 10^{-10}$. In order to illustrate the effect of $\epsilon$, we consider three values of violation $\epsilon \in \{0.2, 0.1, 0.05\}$.
According to Lemma \eqref{lem:epsBetaResult} we require $M \in \{683, 1403, 2842\}$ training samples, respectively, to satisfy the specified violation $\epsilon$ with the desired confidence.
In addition, we consider an unconstrained synthesis problem with $M \in \{683, 1403, 2842\}$ that does not enforce the descent condition \eqref{eq:Scenario_Descent}.
We take $M$ random samples from $\mathcal{D}_{\text{data}}$ and use the learned terminal cost function in our original MPC problem from \eqref{eq:MPC} with a horizon of $N=1$:
\vspace{4pt}
\begin{subequations}
\begin{align}
\hspace{-10mm}
V(x_t) :=& \min_{\mathbf{u}} \, \ell\left(x_{k \mid t}, u_{k \mid t}\right) + \mathcal{V}(x_{k+1 \mid t}; \theta^*_M),\\
\text{s.t. } x_{k+1 \mid t} = & f(x_{k \mid t}, u_{k \mid t} ), \quad \hspace{-3mm} \forall k \in \{0,1, \dots, N - 1\}, \\
u_{k \mid t} \in & \, \mathbb{U}, \quad \hspace{12.5mm} \forall k \in \{0,1, \dots, N - 1\}, \\
x_{k \mid t} \in & \, \mathbb{X}, \quad \hspace{12.5mm} \forall k \in \{1,2, \dots, N\}, \\
x_{0 \mid t} = \, & x_t.
\end{align}
\end{subequations}

\begin{figure*}[t!]
    \centering
    \begin{subfigure}{0.32\textwidth} % Adjust width to fit three in a row
        \centering
        \includegraphics[width=\linewidth]{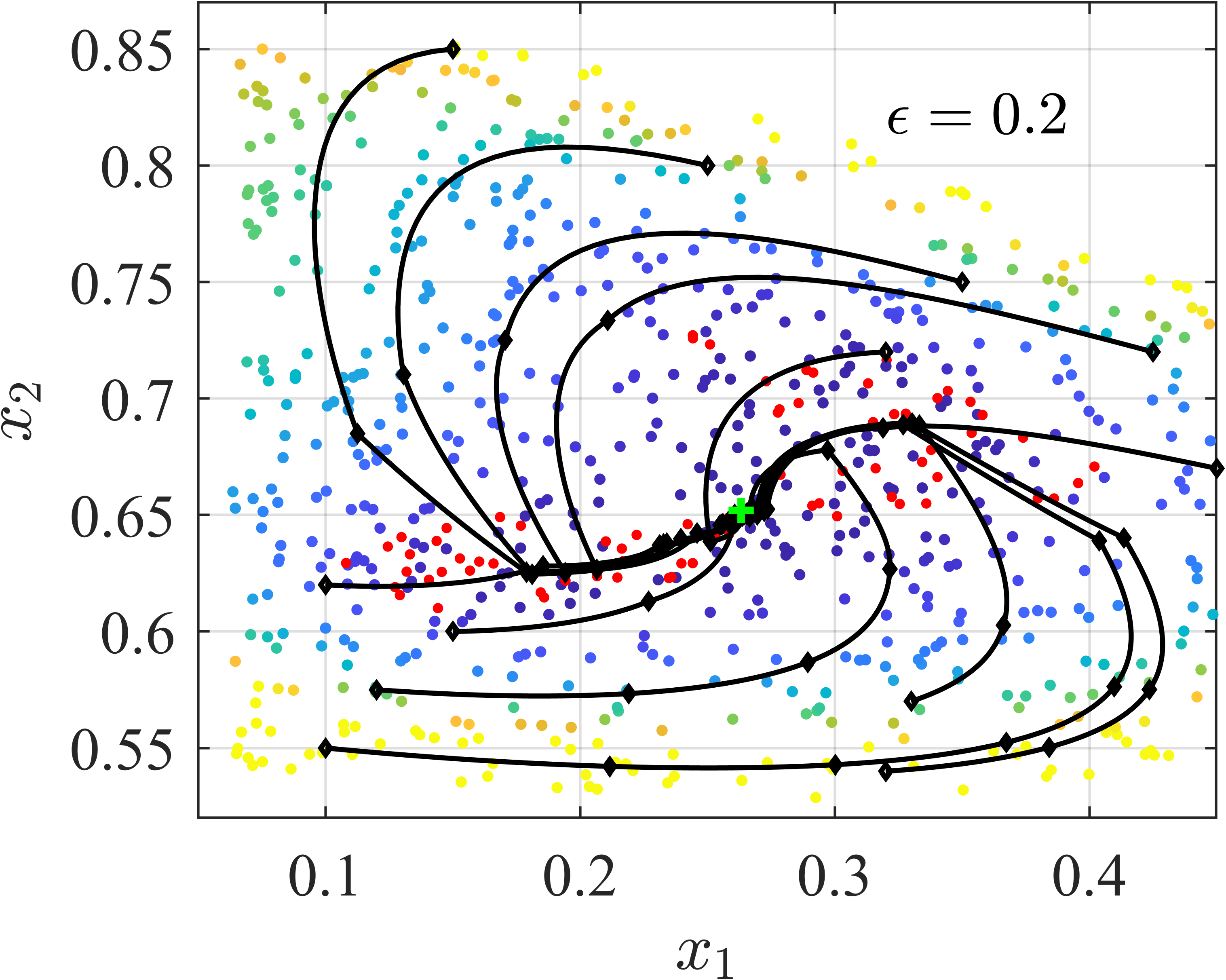}
        \caption{Training points $M=683$,
        Empirical violation $\Tilde{\epsilon} = 0.1091$.}
    \end{subfigure}
    \hfill
    \begin{subfigure}{0.32\textwidth}
        \centering
        \includegraphics[width=\linewidth]{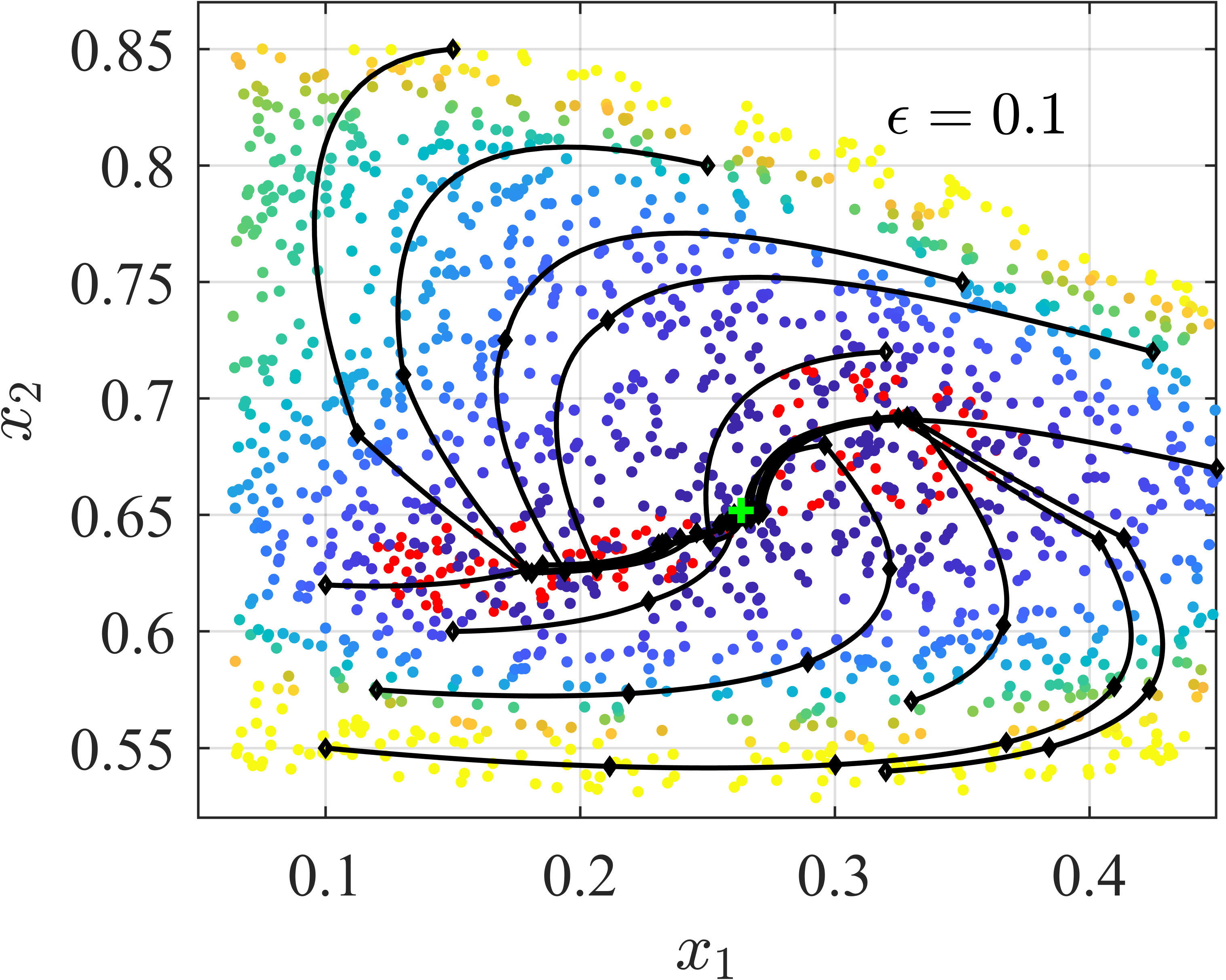}
        \caption{Training points $M=1403$,
        Empirical violation $\Tilde{\epsilon} = 0.0869$.}
    \end{subfigure}
    \hfill
    \begin{subfigure}{0.32\textwidth}
        \centering
        \includegraphics[width=\linewidth]{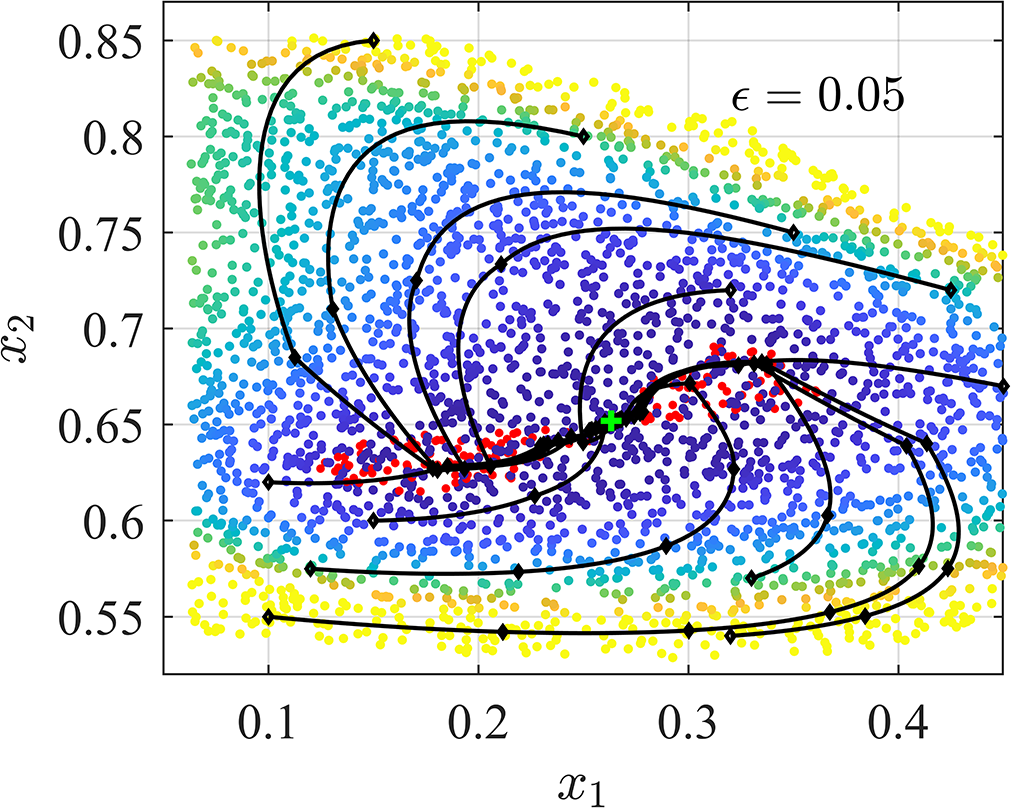}
        \caption{Training points $M=2842$,
        Empirical violation $\Tilde{\epsilon} = 0.0483$.}
    \end{subfigure}
    \caption{Phase plots of twelve closed-loop trajectories of the proposed MPC $(N=1)$ with the learned cost-to-go from $S_M$. The color-graded dots depict the learned cost-to-go $\mathcal{V}(x, \theta^*_M)$ evaluated at the training points. Red dots indicate the violation of the descent condition in the test data set. The green cross indicates the set point of the MPC.}
    \label{fig:viol}
    \vspace{-1mm}
\end{figure*}

\subsection{Numerical Results}
We empirically validate the result of Theorem $\ref{th:Stab}$ by solving a single iteration of the MPC \eqref{eq:MPC} at the $M_{\text{test}}=6226$ states of the test data set $\mathcal{D}_{\text{test}}$ and verifying the basic stability condition \eqref{eq:descentcond}.
In addition, we validate the closed-loop performance of the proposed MPC \eqref{eq:MPC} with the learned terminal cost function $\mathcal{V}(\cdot; \theta^*_M)$ by comparing it against the baseline expert MPC \eqref{eq:Expert_MPC}, shown in Fig. \ref{fig:cstr_full_nmpc}, for various initial conditions.
Computations are performed on an Intel$^{\text{\textregistered}}$ Core$^{\text{\texttrademark }}$ i7 CPU with 32 GB RAM, in MATLAB
using CasADi \cite{andersson_casadi_2019} with the IPOPT \cite{wachter_implementation_2006} solver. The synthesis problem, which is solved offline, was solved in the order of seconds. Expert demonstrations could typically be generated offline, for example, using open-loop MPC solutions, akin to the data generation framework used in direct policy approximation \cite{Allgower_Policy_2018, Krishnamoorthy_Policy_2021}, or from human expert demonstrations.

\begin{table}[b]
    \centering
    \vspace{-4mm}
    \caption{Mean Squared Error $(\mathcal{G})$ of the learned terminal cost function $\mathcal{V}(\theta^*_M)$ and the true cost-to-go $\mathcal{J}$.}
    \begin{tabular}{lccc}
        \toprule
        \textbf{No. training points $(M)$} & \textbf{683} & \textbf{1403} & \textbf{2842} \\
        \midrule
        Constrained $S_M$ & 0.3510 & 0.6047 & 1.2986 \\
        Unconstrained $S_M$   & 0.3134 & 0.5372 & 1.1097 \\
        \bottomrule
    \end{tabular}
    \label{tab:mse}
\end{table}

The unconstrained terminal cost function has in all three cases a lower mean squared error (MSE) with respect to the true cost-to-go $\mathcal{J}$ compared to the constrained terminal cost function, cf. Table \ref{tab:mse}. However, the unconstrained terminal cost function converges to an undesired equilibrium, as seen in Fig. \ref{fig:cstr_nodesc_viol}. 
As mentioned before, simply minimizing the point-wise error to the cost-to-go does not necessarily yield a `good' terminal cost function. Rather the descent along closed-loop trajectories, as enforced in the synthesis program $S_M$, is essential. Consequently, the proposed MPC converges to the desired setpoint, as seen in Fig. \ref{fig:cstr_3d} and \ref{fig:viol}. Note that the constrained radial basis network learns a very complex and nonconvex terminal cost function, Fig. \ref{fig:cstr_3d}.

\begin{table}[b]
    \centering
    \caption{Comparison of Expert and Proposed MPC.}
    \begin{tabular}{lcccc}
        \toprule
        & \multicolumn{1}{c}{\textbf{Expert MPC}} & \multicolumn{3}{c}{\textbf{Proposed MPC}} \\
        \cmidrule(lr){3-5}
        Specified violation $(\epsilon)$ & & $0.2$ & $0.1$ & $0.05$ \\
        \midrule
        Empirical violation $(\Bar{\epsilon})$ & - & 0.109 & 0.087 & 0.048 \\
        No. training points $(M)$ & - & 683 & 1403 & 2842 \\
        Avg. solve time [ms] & 139.5 & 16.6 & 17.2 & 15.0 \\
        Max. solve time [ms] & 570.6 & 52.7 & 100.2 & 47.6 \\
        \bottomrule
    \end{tabular}
    \label{table:results}
\end{table}

The terminal cost function empirically satisfies the basic stability condition everywhere but an $\epsilon$ region, cf. Table \ref{table:results}. Note that with confidence $1-\beta$ it is unlikely that we do not attain this result.
According to Theorem \ref{th:Stab}, the violation set reduces with increasing number of training points, see Fig. \ref{fig:viol}.
For this particular system, the true cost-to-go is relatively flat near the setpoint, as seen in Fig. \ref{fig:cstr_full_nmpc}. Consequently, it is more difficult to enforce the descent condition in this region, leading to more concentrated violations.
Upon close inspection, it can be observed that the violation set in this case is a disjoint set, as there are many points around the origin that are constrained to satisfy the stability condition.
The proposed MPC attains very similar performance to the expert MPC with a significant reduction in the horizon length and hence in computation time, see Table \ref{table:results}.

\subsection{Adaptations of the Optimal Control Problem}
\begin{figure}[b]
    \centering
    \vspace{-4.5mm}
    \includegraphics[width=0.5\textwidth]{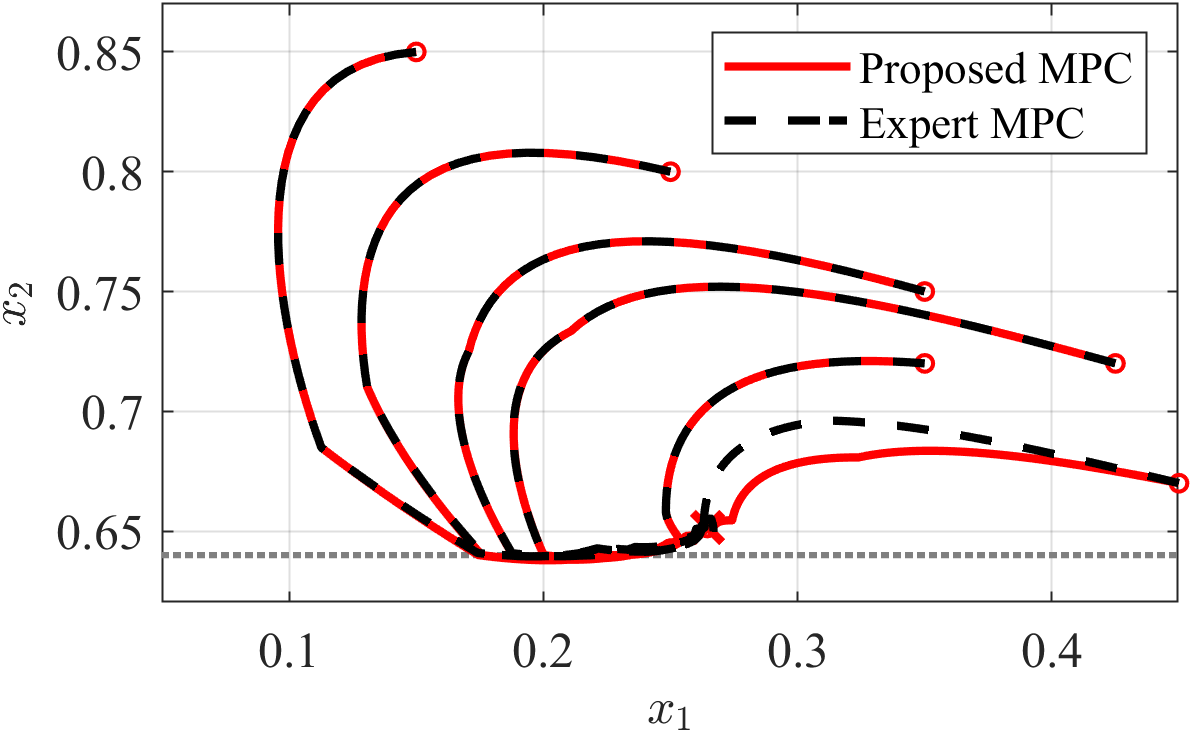}
    \caption{Phase plot of the proposed MPC $(N=1)$ trained with \text{$M = 2842$} training points and the expert MPC $(T=50)$. Soft-constraints are used for the adapted state constraint which is indicated by the dotted line.}
    \label{fig:nmpc_adapted}
\end{figure}

As discussed in Section \ref{sec:Intro}, a strong advantage of implicit MPC over explicit MPC is its flexibility and robustness against modifications in the optimal control problem. To this end, we demonstrate the proposed MPC in a modified CSTR problem with \textit{tightened} state constraints, namely, we tighten the lower bound of $x_2$ from $0.4519$ to $0.64$. In this example we use a learned terminal cost function with $M=2842$ without modifying or re-training the cost function for the adapted state constraints, i.e., we only adjust the state constraint set $\mathbb{X}$. We use soft-constraints on the states to retain recursive feasibility of the MPC problem. 
It is important to note that the tightening of the constraints may increase the probability of entering the violation set.
Figure \ref{fig:nmpc_adapted} shows the results of the proposed MPC $(N=1)$ compared against the expert MPC $(T=50)$, both having tightened soft-constraints. The proposed MPC attains almost identical performance to the expert MPC with only minor violations of the state constraints, demonstrating the flexibility and performance of implicit value function-augmented MPC.
\vspace{4.5mm}

\section{Discussion}
In this work, we presented a data-driven synthesis method of a terminal cost function $\mathcal{V}$ for a general nonlinear MPC.
If we can solve the synthesis problem, then we can have confidence that the resulting function is an appropriate terminal cost function.
However, the probabilistic certificates obtained by the scenario program only hold under the same distribution over $\mathcal{X}$, namely uniformly, and hence have some theoretical limitations.
For example, under the closed loop policy $\kappa_N$, the state space will not be sampled uniformly. As such, the violation probability $\epsilon$ does not necessarily carry over to the decrease of the value function $V(x_t)$ when applied in closed-loop.
In future work, we aim to account for these distribution shifts in our analysis. Furthermore, the violation set $\Delta_\epsilon$ could steer the system to a (local) minimum of $\mathcal{V}$.
However, note that the results presented in this manuscript only rely on the structure of the scenario program, namely convexity, and the optimality of the MPC. We aim to address the closed-loop analysis with the learned terminal cost in future work by exploiting the properties of $\mathcal{V}$, the dynamics $f$ and the MPC policy $\kappa_N$.
As mentioned in Section \ref{sec:Intro}, this work is a first step in this direction and our aim is to address the current limitations in future research.

\section{Conclusions}
This paper presents a novel method to synthesize a terminal cost function for nonlinear MPC using supervised learning. The method uses (expert) demonstrations in the form of state-input-cost tuples to solve a convex synthesis problem. We use scenario optimization to guarantee that the \textit{basic stability condition} of the terminal cost function is satisfied almost everywhere with high probability.
We demonstrated the methodology in a numerical example where the cost-to-go is approximated with a radial basis network. 

Firstly, the learned terminal cost function enables a reduction in the horizon to a one-step MPC with minimal performance loss, and is flexible to modifications of the OCP.
Secondly, in the context of using MPC as a local approximation of a globally optimal control policy, this methodology enables the use of more complex and expressive MPC value functions, while enforcing stabilizing properties.
The next steps for value function-augmented MPC include
(i) analysis of the closed-loop MPC including recursive feasibility,
\text{(ii) handling} of distribution shifts in training data, and
\mbox{(iii) robustification} against disturbances and modeling errors.

\section*{Acknowledgements}
The authors would like to thank Prof. Simone Garatti from Politecnico di Milano for the inspiring discussions.

\newpage
\bibliographystyle{IEEEtran}
\bibliography{Ref}

\end{document}